\newtheorem{myth}{Theorem}
\newtheorem{myle}{Lemma}
\newtheorem{myprop}{Proposition}

\documentclass[a4paper,10pt,greek,english]{article}
\usepackage{babel}
\usepackage[iso-8859-7]{inputenc}
\newfont{\Bdd}{msbm8 scaled\magstep1}
\begin{document}
\title{\bf An Effective Version of Chevalley-Weil Theorem  for Projective Plane Curves}
\author{Konstantinos Draziotis   and Dimitrios Poulakis }
\date{}
\maketitle

\begin{abstract}
We obtain a quantitative version of the classical Chevalley-Weil
theorem for curves. Let $\phi : \tilde{C} \rightarrow C$ be an
unramified morphism of non-singular plane projective curves
defined over a number field $K$. We calculate an effective upper
bound for the norm of the relative discriminant of the number
field $K(Q)$ over $K$ for any point $P\in C(K)$ and
$Q\in{\phi}^{-1}(P)$.

\ \\
2000 MCS: 14G25, 14H25, 11G30.
\end{abstract}

\section{Introduction}
Let $\phi : V \rightarrow W$ be an unramified covering of projective normal varieties defined over
a number field $K$. By the classical theorem of Chevalley-Weil  \cite{Chevalley}, \cite{Weil},
\cite[Theorem 8.1, page 45]{Lang 83},  \cite[page 292]{Hindry}, there exists a finite extension
$L/K$ such that $\phi^{-1}(W(K)) \subseteq V(L)$.
In \cite[Theorem 1.1]{Draziotis}, we obtained  a quantitative version of the
 Chevalley-Weil theorem in case where
 $\phi : \tilde{C} \rightarrow C$ is an unramified morphism
  of non-singular affine plane curves defined over $K$.
More precisely, we gave, following a new  approach, an effective
upper bound for the  relative discriminant of the minimal field of
definition $K(Q)$ of $Q$ over $K$  for any integral point $P\in
C(K)$ and $Q \in{\phi}^{-1}(P)$. In this paper, we  consider the
case where $\phi : \tilde{C} \rightarrow C$ is an unramified
morphism of non-singular projective plane curves defined over $K$
and we obtain, extending our method,  an effective upper bound for
the relative discriminant of  $K(Q)$ over $K$  for any $P\in C(K)$
and $Q\in{\phi}^{-1}(P)$.

Consider the set of absolute values on ${\mbox{\Bdd Q}}$
consisting of the ordinary absolute value and for every prime $p$
the $p$-adic absolute value $|\cdot|_p$ with $|p|_p  = p^{-1}$.
Let $M(K)$ be a set of symbols $v$ such that with every $v \in
M(K)$ there is precisely one associated absolute value
$|\cdot|_{v}$ on $K$ which extends one of the above absolute
values of $\mbox{\Bdd Q}$. We denote by $d_v$  its local degree.
Let ${\bf x} = (x_0:\ldots:x_n)$ be  a point of the projective
space ${\mbox{\Bdd P}}^{n}(K)$ over $K$.  We define the field
height $H_{K}(\bf x)$ of $\bf x$ by
$$ H_{K}({\bf x}) = \prod_{v \in M(K)} \max\{|x_{0}|_{v},\ldots,|x_{n}|_{v}\}^{d_v}.$$
Let   $d$ be the degree of $Κ$. We define  the absolute height
$H(\bf x)$ by $H({\bf x}) = H_{K}({\bf x})^{1/d}$. Furthermore,
for $x \in K$ we put $H_K(x) = H_K(1:x)$ and $H(x) = H(1:x)$. If
$G \in K[X_1,\ldots,X_m]$, then we define the field height
$H_K(G)$ and the absolute height $H(G)$ of $G$ as the field height
and the absolute height of the point whose coordinates are the
coefficients of $G$. For an account of the properties of heights
see \cite[chapter VIII]{Silverman} or \cite[chapter 3]{Lang 83}.

Let $\overline{K}$ be an algebraic closure of $K$ and $O_K$ the
ring of algebraic integers of $K$. If $M$ is a finite extension of
$Κ$, then we denote by $D_{M/K}$ the relative discriminant of the
extension $M/K$ and  by  $N_M$ the norm from  $M$ to ${\mbox{\Bdd
Q}}$.

Let $F,\bar{F} \in K[X_1,X_2,X_3]$ be two homogeneous absolute
irreducible polynomials with $N = \deg F > 1$ and $\bar{N} =  \deg
\bar{F} > 1$. We denote by $C$ and $\bar{C}$  the projective
curves defined by  $F(X_1,X_2,X_3) = 0$ and $\bar{F}(X_1,X_2,X_3)
= 0$ respectively. Let $\phi : \bar{C} \rightarrow C$ be a
nonconstant morphism of degree $m > 1$ defined by
$\phi(X_1,X_2,X_3) =
({\phi}_1(X_1,X_2,X_3),{\phi}_2(X_1,X_2,X_3),{\phi}_3(X_1,X_2,X_3))$,
where ${\phi}_i(X_1,X_2,X_3)$ $(i=1,2,3)$ are relatively prime
homogeneous polynomials in  $K[X_1,X_2,X_3]$  of the same degree
$M$. Let  $\Phi$ be a point in the projective space having as
coordinates the coefficients of ${\phi}_i$ $(i = 1,2,3)$.

\begin{myth}
 Suppose that $C$ is nonsingular
 and the morphism $\phi : \bar{C} \rightarrow C$  unramified.  Then for any
point $P \in C(K)$  and $Q \in {\phi}^{-1}(P)$, we have
$$N_K(D_{K(Q)/K}) <    \Omega (H(F)^{6N^2\bar{N}} H(\Phi_i)^{\bar{N}}
H(\bar{F})^M)^{\omega d m^3 M^7N^{30}\bar{N}^{13}},$$ where
$\Omega$ is an effectively computable constant in terms of
 $N,\bar{N},M,m$ and $d$, and $\omega$ a numerical constant.
\end{myth}
{\it Remarks.}
1) By \cite[Corollary 3, p. 120]{Shafarevich}, the curve $\bar{C}$  is nonsingular. \\
2) Since $m > 1$, the quantity $M$  is $> 1$. \\
3) By Hurwitz' s formula,   $\bar{C}$ and $C$  have positive genus and  $\bar{N} \geq N \geq 3$. \\
4) Since $\bar{F}(X,Y,Z)$ divides
$F(\phi_1(X,Y,Z),\phi_2(X,Y,Z),\phi_3(X,Y,Z))$,  $H(\bar{F})$ and
$\bar{N}$ can be bounded  by  constants depending only on  $F$ and
$\phi$.

Let $K(\bar{C})$ and $K(C)$ be the function fields of $\bar{C}$
and $C$, respectively, over $K$, $P = (p_1:p_2:p_3)$ and $\phi^{*}
: K(C) \rightarrow K(\bar{C})$ the field homomorphism associated
to $\phi$. We denote by $\phi_{j,i}$ the function on $\bar{C}$
defined by the fraction $\phi_j/\phi_i$. The idea of the proof of
Theorem 1 is as follows. For every affine view $C_i$, with $X_i =
1$ $(i = 1,2,3)$, of $C$ we construct two primitive elements
$u_{is}$ $(s=1,2)$ for the field extension
$K(\bar{C})/\phi^{*}(K(C))$ which are  integral
 over the ring $K[\phi_{j,i},\phi_{k,i}]$ and such that
$K(u_{is}(Q)) = K(Q)$. Further, we construct polynomials
$P_{is}(X,Y,U)$ $(s=1,2)$  representing the minimal polynomials of
$u_{is}$ over $K[\phi_{j,i},\phi_{k,i}]$ such that the
discriminants $D_{is}(X,Y)$ of $P_{is}(X,Y,U)$ $(s = 1,2)$ have no
common zero on $C_i$. It follows that for every prime ideal $\wp$
of $O_K$ with quite large norm there is $i \in \{1,2,3\}$ such
that $\wp$ cannot divide both $D_{is}(p_j/p_i,p_k/p_i)$ $(s=1,2)$
and hence cannot divide the discriminant of $K(Q)$. Thus, we
determine the prime ideals of $K$ which are ramified in $K(Q)$ and
the result follows. A totally different effective approach of
Chevalley-Weil theorem  is given in \cite[Chapter 4]{Bilu}.

The paper is organized as follows.  In section 2 we give some
auxiliary results and in section 3 we obtain the proof of Theorem
1.

{\it Notations.}  If $C$ is a projective plane curve defined over
$\overline{K}$, then we denote by $O(U)$ the ring of regular
functions on an open subset $U$ of $C$ and by $\overline{K}(C)$
the function field of $C$. Let $G$ be a homogeneous polynomial of
$\overline{K}[X_1,X_2,X_3]$. We denote by $D_C(G)$ and $V_C(G)$
the set of points $P \in C(\overline{K})$ with $G(P) \neq 0$ and
$G(P) = 0$ respectively. Finally, throughout the paper, we denote
by $\Lambda_1(a_1,\ldots,a_s)$, $\Lambda_2(a_1,\ldots,a_s), \ldots
$ effectively computable positive numbers in terms of indicated
parameters.

\section{Auxiliary Results}

We keep the notations and the assertions of the Introduction.
 The restriction of $\phi$ on  ${\phi}^{-1}(D_C(X_i))$ is a finite
morphism. Thus, the  associated ring homomorphism ${\phi}^{*} :
O(D_C(X_i)) \rightarrow O({\phi}^{-1}(D_C(X_i)))$, defined by
$\phi^{*}(f) = f{\circ}{\phi}$, for every $f \in O(D_C(X_i))$,  is
injective and the ring $O({\phi}^{-1}(D_C(X_i)))$ is finite over
$\phi^{*}(O(D_C(X_i))$. We denote by $\bar{x}_{j,i}$ and $x_{j,i}$
the functions defined by $X_j/X_i$ on $\bar{C}$ and $C$,
respectively. The function ${\phi}^{*}(x_{j,i})$ is defined by the
fraction $\phi_j/\phi_i$ and so  ${\phi}_{j,i} =
{\phi}^{*}(x_{j,i})$. Then we have $\phi^{*}(O(D_C(X_i)) =
 \overline{K}[{\phi}_{j,i},{\phi}_{k,i}]$. Let $\rho$ be an integer
such that for every $(z_1:z_2:z_3) \in  V_{\bar{C}}(X_i)$ we have
$z_k +\rho z_j \neq 0$, where $\{i,j,k\} = \{1,2,3\}$ with $j <
k$. Thus, the poles of the function $u =
\bar{x}_{k,i}+\rho\bar{x}_{j,i}$ are the points of
$V_{\bar{C}}(X_i)$. Put $\Pi_i = {\phi}^{-1}(D_C(X_i))\cap
V_{\bar{C}}(X_i)$.

\begin{myprop}
There is a monic polynomial $f(T) \in K[T]$ such that the function
$\tilde{u} = uf({\phi}_{j,i})$ is integral over
$K[{\phi}_{j,i},{\phi}_{k,i}]$. We have $\deg f \leq \bar{N}$,
$$H(f) < \Lambda_1(\rho,M,N,\bar{N}) H(F)^{\bar{N}} H(\bar{F})^{MN}
H(\Phi)^{N\bar{N}},$$
 and the roots of $f(T)$ are
the elements $\phi_{j,i}(R)$, where $R \in \Pi_i$. Moreover, there
is a polynomial of $K[X_j,X_k]$,
$$P(X_j,X_k,U) = U^{\mu}+p_1(X_j,X_k)U^{\mu-1}+{\cdots}+p_{\mu}(X_j,X_k),$$
such that $P({\phi}_{j,i},{\phi}_{k,i},U)$ is the minimal
polynomial of $\tilde{u}$ over $K[{\phi}_{j,i},{\phi}_{k,i}]$. We
have $\mu  \leq m$, $\deg p_l < 11MN^4\bar{N}^2$ $(l =
1,\ldots,\mu)$ and
$$H(P) < {\Lambda}_2(\rho,m,M,N, \bar{N})
(H(F)^{6N^2\bar{N}}H(\Phi)^{\bar{N}}
H(\bar{F})^{M})^{240mM^3N^{12}\bar{N}^5}.$$
\end{myprop}

For the proof of Proposition 1 we shall need the following lemma.

\begin{myle}
 There is a polynomial  $G(W,X,U) \in K[W,X,U]\setminus \{0\}$ such that
$G({\rho},\phi_{j,i},u) = 0$. We have $\deg_X G \leq N\bar{N}$,
$\deg_U G \leq 2MN \bar{N}$,  $\deg_W G \leq 2MN \bar{N}$ and the
polynomial $G_{\rho}(X,U) = G(\rho,X,U)$ satisfies
$$H(G_{\rho}) <
\Lambda_3(\rho,M,N,\bar{N}) H(F)^{\bar{N}} H(\bar{F})^{MN}
H(\Phi)^{N\bar{N}}.$$
\end{myle}
{\it Proof.} We may suppose, without loss of generality, that $j =
1$, $k = 2$ and $i = 3$. Consider the polynomials
$\bar{F}_1(W,V,U) = \bar{F}(V,U-WV,1)$ and
$$E(W,X,V,U) =
F(X\phi_3(V,U-WV,1),\phi_2(V,U-WV,1),\phi_3(V,U-WV,1)).$$
 We have $\bar{F}_1(\rho,\bar{x}_{1,3},u) =
E(\rho,\phi_{1,3},\bar{x}_{1,3},u_{\rho}) = 0.$  If $G(W,X,U)$ is
the resultant of $E(W,X,V,U)$ and $\bar{F}_1(W,V,U)$ with respect
to $V$, then  $G({\rho},{\phi}_{1,3},u) = 0$.

Suppose that $G(W,X,U)$ is equal to zero. Thus, since
$\bar{F}_1(W,V,U)$ is absolutely irreducible,  $\bar{F}_1(W,V,U)$
divides $E(W,X,V,U)$. It follows that $\bar{F}(V,U,1)$ divides
$F(X\phi_3(V,U,1),\phi_2(V,U,1),\phi_3(V,U,1))$. Write
$$F(X_1,X_2,X_3) = A_0(X_2,X_3)X_1^n+\cdots +
A_n(X_2,X_3),$$ where $A_i(X_2,X_3)$ $(i=0,\ldots,n)$ are
homogeneous polynomials with $\deg A_i  = N-n+i$. If $P =
(p_1:p_2:1) \in D_{\bar{C}}(\phi_3)$, then
$$A_0(\phi_{2,3}(P),1)(X_1/\phi_3(P))^n+\cdots +
A_n(\phi_{2,3}(P),1) = 0.$$ It follows that $A_j(\phi_{2,3}(P),1)
= 0$ $(j = 0,\ldots,n)$
 which is a contradiction since $F(X_1,X_2,X_3)$ is absolutely irreducible.
 Thus $G(W,X,U)$ is not zero.

By  \cite[Lemma 4.2]{Draziotis}, we have $\deg_X G \leq N
\bar{N}$, $\deg_U G \leq 2MN \bar{N}$, and $\deg_W G \leq 2MN
\bar{N}$. Further, if $G_{\rho}(X,U) = G({\rho},X,U)$,
$E_{\rho}(X,V,U) = E({\rho},X,V,U)$ and $\bar{F}_{\rho}(V,U) =
\bar{F}_1(\rho,V,U)$, then
$$H(G_{\rho}) < \Lambda_4(M,N,\bar{N}) H(E_{\rho})^{\bar{N}}
H(\bar{F}_{\rho})^{MN}.$$

By  \cite[Lemma 4.4]{Draziotis}, we obtain
$$H(\bar{F}_{\rho}) \leq
2^{\bar{N}} ({\bar{N}}+1)\  {\max}\{1,|{\rho}|\}^{\bar{N}}
H({\bar{F}}).$$ Next, put $\varphi_{\rho,l}(V,U) = \phi_l(V,U-\rho
V)$ $(l=1,2)$. By \cite[Lemma B.7.4]{Hindry}, for every absolute
value $|\cdot|_v$ of $K$,
$$|E_{\rho}|_v \leq \max\{1,|2N|_v^2\} |F|_v \max_{0\leq j \leq N}\{|\varphi_{\rho,2}^j|_v|\varphi_{\rho,3}^{N-j}|_v\}$$
 and for every positive number $k$,
$$|\varphi_{\rho,l}^k|_v \leq  \max\{1,|2M|_v\}^{2(k-1)M}
|\varphi_{\rho,l}|_v^k.$$ Furthermore, the proof of \cite[Lemma
4.4]{Draziotis} gives
$$|\varphi_{\rho,l}|_v \leq \max\{1,|\rho|_v\}^M \max\{1,|2|_v\}^M
\max\{1,|M+1|_v\} |\phi_l|_v  \ \ (l=1,2).$$ The above
inequalities yield
$$H(E_{\rho}) < \Lambda_5(\rho,M,N,\bar{N}) H(F) H(\Phi)^N.$$
Combining all theses estimates, the bound for $H(G_{\rho})$
follows.

\vspace{3mm} {\it Proof of Proposition 1.} By Lemma 1, there is
$G_{\rho}(X,U) \in K[X,U]$  such that $G_{\rho}(\phi_{j,i},u) =
0$. Write $G_{\rho}(X,U) = g_0(X)U^{\nu}+\cdots+g_{\nu}(X)$. Thus,
$ug_0({\phi}_{j,i})$ is an integral element over
$K[{\phi}_{j,i},{\phi}_{k,i}]$ and so $ug_0({\phi}_{j,i}) \in
O({\phi}^{-1}(D_C(X_i)))$.

If $h \in \overline{K}(\bar{C})$ and $S \in \bar{C}$, then we
denote by ${\rm ord}_S(h)$ the order of  $h$ at $S$. Put $B_R =
\phi_{j,i}(R)$, where $R \in \Pi_i$. Let $m_R$ be the smallest
integer such that $(\phi_{j,i}-B_R)^{m_R}u$ is defined at $R$.
Then $m_R \leq |{\rm ord}_R(u)|$. Set $f(X) = \prod_{R\in
\Pi_i}(X-B_R)^{m_R}$. We have $uf(\phi_{j,i}) \in
O({\phi}^{-1}(D_C(X_i))$ and since
$[\overline{K}(\bar{C}):\overline{K}(u)] = \bar{N}$, we obtain
$\deg f = \sum_{R\in \Pi_i} m_R \leq \bar{N}$.  The elements of
the Galois group ${\rm Gal}(\overline{K}/K)$ permute the elements
of $\Pi_i$ and consequently the numbers $B_R$. For every $\sigma
\in {\rm Gal}(\overline{K}/K)$, we have ${\rm
ord}_R(\phi_{j,i}-B_R)= {\rm
ord}_{R^{\sigma}}(\phi_{j,i}-B_{R^{\sigma}})$ and ${\rm ord}_R(u)=
{\rm ord}_{R^{\sigma}}(u)$. It follows that   $m_R =
m_{R^{\sigma}}$. Hence $f(X) \in K[X]$. Since $ug_0({\phi}_{j,i})
\in O({\phi}^{-1}(D_C(X_i)))$, we have  $g_0(X) = f(X)l(X)$, where
$l(X) \in K[X]$. By \cite[Proposition B.7.3]{Hindry},
 $H(f) \leq e^{N \bar{N}} H(G_{\rho})$. The bound for $H(f)$ follows.

Consider the polynomial
$$\tilde{G}_{\rho}(X,U) = l(X)U^{\nu}+g_1(X)U^{\nu-1}+g_2(X)f(X)U^{\nu-1}+\cdots+g_{\nu}(X)f(X)^{\nu-1}.$$
We have $\tilde{G}_{\rho}({\phi}_{j,i},uf({\phi}_{j,i}) = 0$. The
estimates for $G_{\rho}(X,U)$ and \cite[Proposition B.7.4]{Hindry}
yield
$$H(\tilde{G}_{\rho}) < \Lambda_7(\rho,M,N,\bar{N}) (H(F)^{\bar{N}} H(\bar{F})^{MN}
H(\Phi)^{N\bar{N}})^{2MN\bar{N}}.$$ Using \cite[Proposition
2.1]{Draziotis} and the estimates for $\tilde{G}_{\rho}$,  we
obtain the existence of polynomial $P(X_j,X_k,U) \in K[X_j,X_k,U]$
having the required properties.

\begin{myle}
 Let $P{\in}C(K)$ and $Q \in {\bar{C}}(\overline{K})$ with $\phi (Q) = P$. Then
$$N_K(D_{K(Q)/K}) < ((e^3(M+\bar{N}))^{dM\bar{N}} (H_K(P)H_K({\Phi}))^{\bar{N}}
H_K({\bar{F}})^{M})^{40dM^3\bar{N}^3}.$$
\end{myle}
{\it Proof.} We may suppose, without loss of generality, that $Q =
(q_1:q_2:1)$ and $P = (p_1:p_2:1)$. Put $G_1(X_1,U,V) =
X_1{\phi}_3(U,V,1)-{\phi}_1(U,V,1).$ Then $G_1(p_1,q_1,q_2) =
\bar{F}(q_1,q_2,1) = 0.$ We denote by $R_1(U)$ and $R_2(V)$  the
resultants of $\bar{F}(U,V,1)$ and ${\Gamma}(U,V) = G_1(p_1,U,V)$
with respect to $V$ and $U$. Then $R_1(q_1) = R_2(q_2) = 0$. By
\cite[Lemma 4.2]{Draziotis} and \cite[Proposition
B.7.4(b)]{Hindry}  we obtain
$$H(R_i) \leq (M+\bar{N})! (\bar{N}+1)^{M} (M+1)^{\bar{N}} (2H(p_{1})H({\Phi}))^{\bar{N}} H({\bar{F}})^{M}.$$
Furthermore, we have $\deg R_i \leq 2M\bar{N}$.

Let $B_i(T) = T^{m_i}+b_1T^{m_i-1}+ \cdots +b_{m_i}$, where $m_i
\leq 2M\bar{N}$,  be the  irreducible polynomial of $q_i$ over
$K$.  By \cite[Lemma 4.1]{Draziotis} there is a positive integer
$\beta_i$ with $\beta_i \leq H_K(B_i)^{m_i}$ such that $\beta_i
b_1 \cdots b_{m_i} \in O_K$. Then $\beta_iq_i$ is an algebraic
integer with minimal polynomial $\bar{B}_i(T) = T^{m_i}+\beta
b_1T^{m_i-1}+ \cdots +{\beta}^{m_i} b_{m_i}$. Using
\cite[Proposition B.7.3]{Hindry} we obtain
$$H(\bar{B}_i) \leq H(B_i) \beta_i^{m_i} \leq (e^{2M\bar{N}}
H(R_i))^{1+2dM\bar{N}}. $$
 Let $\Delta(\bar{B}_i)$ be the
discriminant of $\bar{B}_i(T)$. By \cite[Lemma 5]{Poulakis1}, we
have
$$N_K(\Delta(\bar{B}_i)) \leq H_K(\Delta(\bar{B}_i)) \leq   m_i^{3m_id} H_K(\bar{B}_i)^{2m_i-2} \leq (e^{2dM\bar{N}}H_K(R_i))^{9dM^2\bar{N}^2}.$$

Put $K_i = K(q_i)$. Since $b_iq_i$ is an algebraic integer, the
discriminant $D_i$ of the extension $K_i/K$ divides the
discriminant of  $1,b_iq_i,{\ldots},(b_iq_i)^{m_i-1}$ which is
equal to  $\Delta(\bar{B}_i)$. Thus $N_K(D_i) \leq
|N_K(\Delta(\bar{B}_i))|$. If $I(T)$ is the irreducible polynomial
of $b_2q_2$ over $K_1$, then $I(T)$ divides $\bar{B}_2(T)$ (in
$K_1[T]$) and so the discriminant ${\Delta}(I)$ of $I(T)$ divides
$\Delta(\bar{B}_2)$. Hence, $D_{K(Q)/K_{1}}$  divides
$\Delta(\bar{B}_2)$. Thus,
$$N_K(D_{K(Q)/K}) \leq N_K(D_1)^{2M\bar{N}} N_{K_1}(D_{K(Q)/K_{1}})
\leq (N_K(\Delta(\bar{B}_1) N_K(\Delta(\bar{B}_2))^{2M\bar{N}}.$$
Using the upper bounds for  $N_K(\Delta(\bar{B}_i))$ and
$H_K(R_i)$, the result follows.

\section{Proof of Theorem 1}
Let $P = (a_1:a_2:a_3)$, $Q \in {\phi}^{-1}(P)$ and  $L = K(Q)$.
 If  $a_{j} = 0$ for some $j{\in}\{1,2,3\}$, then  \cite[Lemma 4]{Poulakis2} gives $H(P) < 2H(F)$. So
 Lemma 2 yields a sharper bound for $N_K(D_{L/K})$ than that
of Theorem 1. Thus, we may suppose that $a_j \neq 0$ $(j =
1,2,3)$.

Let  $\Theta_i$ be the set of $\rho \in \mbox{\Bdd Z}$ such that
for every $(z_1:z_2:z_3) \in V_{\bar{C}}(X_i)$ we have $z_k +\rho
z_j = 0$, where $\{i,j,k\} = \{1,2,3\}$ with $j < k$.  Set
$u_{\rho,i} = \bar{x}_{k,i}+\rho\bar{x}_{j,i}$, where $\rho \not
\in \Theta_i$.
 By Proposition 1, there is a monic polynomial $f_i \in K[T]$ such that the function
$\tilde{u}_{\rho,i} = u_{\rho,i}f_i(\phi_{j,i})$ is integral over
$K[\phi_{j,i},\phi_{k,i}]$, $\deg f_i \leq \bar{N}$, the roots of
$f_i(T)$ are the elements $\phi_{j,i}(R)$, where $R \in
{\phi}^{-1}(D_C(X_i))\cap V_{\bar{C}}(X_i)$ and
$$H(f) < \Lambda_1(\rho,M,N,\bar{N}) H(F)^{\bar{N}} H(\bar{F})^{MN}
H(\Phi)^{N\bar{N}}.$$
 Moreover, there is a polynomial of $K[X_j,X_k,U]$,
$$P_{\rho,i}(X_j,X_k,U) = U^{\mu}+p_{\rho,i,1}(X_j,X_k)U^{\mu-1}+{\cdots}+p_{\rho,i,\mu}(X_j,X_k),$$
such that $P_{\rho,i}(\phi_{j,i},\phi_{k,i},U)$ is the minimal
polynomial of $\tilde{u}_{\rho,i}$ over
$K[\phi_{j,i},\phi_{k,i}]$.  We have $\mu  \leq m$, $\deg
p_{\rho,i,l} <  11MN^4\bar{N}^2$ $(l = 1,\ldots,\mu)$, and
$$H(P) < {\Lambda}_2(\rho,m,M,N, \bar{N})
(H(F)^{6N^2\bar{N}}H(\Phi)^{\bar{N}}
H(\bar{F})^{M})^{240mM^3N^{12}\bar{N}^5}.$$

Suppose that there is $i \in \{1,2,3\}$ such that $f_i(a_j/a_i) =
0$. By  \cite[Lemma 4]{Poulakis2} and \cite[Lemma 7]{Poulakis1},
we have
$$H(P) \leq H(a_j/a_i) H(a_k/a_i) \leq 2(N+1) H(F) (2H(f_i))^{N+1}.$$
Using the bound for  $H(f_i)$, Lemma 2 gives  a sharper bound for
$N_K(D_{L/K})$ than that of Theorem 1. Next, suppose that for
every $i=1,2,3$ we have $f_i(a_j/a_i) \neq 0$ and so $u_{\rho,i}$
is defined at $Q$.

The monomorphism ${\phi}^{*} : O(D_C(X_i)) \rightarrow
O({\phi}^{-1}(D_C(X_i)))$ extends to a field homomorphism
${\phi}^{*} : \overline{K}(C) \rightarrow \overline{K}(\bar{C})$.
We have $\phi^{*}(\overline{K}(C)) =
\overline{K}(\phi_{j,i},\phi_{k,i})$.
 If ${\sigma}_1, \ldots,{\sigma}_m$ are all the
$\overline{K}(C)$-embeddings of $\phi^{*}(\overline{K}(\bar{C}))$
into an algebraic closure of $\phi^{*}(\overline{K}(C))$, then we
denote by $\Gamma_i$ the set of  integers $\rho \not \in \Theta_i$
with $\sigma_{p}(\tilde{u}_{\rho,i}) \neq
\sigma_{q}(\tilde{u}_{\rho,i})$ for $p \neq q$. For every
$\rho{\in}\Gamma_i$, we have $\overline{K}(\bar{C}) =
\phi^{*}(\overline{K}(C))(\tilde{u}_{\rho,i})$ and so $m = \mu$.
Note that at most $m(m-1)/2+ \bar{N}$ integers $\rho$ do not lie
in $\Gamma_i$. Further, there are at most $m(m-1)/2+\bar{N}$
integers $\rho$ such that  $K(u_{\rho,i}(Q)) \neq  K(Q)$. Hence,
there is $r(i)\in {\mbox{\Bdd Z}}$ with $r(i) \in \Gamma_i$ and
$|r(i)| \leq \bar{N}+m^2/2$ such that $K(u_{r(i),i}(Q)) = K(Q)$.

Putting $X_i = 1$ in $F(X_1,X_2,X_3)$ we obtain $F_i(X_j,X_k)$,
with $j <k$. Let $D_{\rho,i}(X_j,X_k)$ be the discriminant  of
$P_{\rho,i}(X_j,X_k,U)$  with respect to $U$. We have $\deg
D_{\rho,i} < 11(2m-1)MN^4\bar{N}^2$. Since
$P_{\rho,i}(\phi_{j,i},\phi_{k,i},U)$ is irreducible, $F_i$ does
not divide  $D_{\rho,i}$. We denote by $J_{r(i),i}$ the set of
points $(z_1:z_2:z_3) \in D_C(X_i)$ with $z_i = 1$ and
$D_{r(i),i}(z_j,z_k) = 0$. By B\'{e}zout's theorem, $|J_{r(i),i}|
< 11(2m-1)MN^5\bar{N}^2$. Thus, if $B_i = J_{r(i),i}
 \cup \{P\}$, then there is an integer $s(i)$ with
$|s(i)| \leq 11m^2\bar{N}^2N^5M$ such that $B_i \cap
\phi(V_{\bar{C}}(X_k+s(i)X_j)) = \O$.

We denote by $\tilde{F}_i(Y_1,Y_2,Y_3)$ and
$\tilde{\phi}_{i,l}(Y_1,Y_2,Y_3)$ the polynomials obtained from
$\bar{F}(X_1,X_2,X_3)$ and $\phi_l(X_1,X_2,X_3)$, respectively,
using the projective change of coordinates $\chi$ defined by $ Y_j
= X_i ,   \   Y_k = X_j,  \       Y_i =  X_k+s(i)X_j$. Set
$\tilde{Q} = \chi(Q)$.  Let $\tilde{C}_i$ be the curve defined by
$\tilde{F}_i(Y_1,Y_2,Y_3) = 0$. The morphism $\psi_i : \tilde{C}_i
\rightarrow C$, defined by $\psi_i(Y_1,Y_2,Y_3) =
(\psi_{i,1}(Y_1,Y_2,Y_3),\psi_{i,2}(Y_1,Y_2,Y_3),\psi_{i,3}(Y_1,Y_2,Y_3))$
is unramified of degree $m$. We denote by $\Psi_i$  a point in the
projective space with coordinates the coefficients of $\psi_{i,s}$
$(s = 1,2,3)$.

Let $y_{j,i}$ be the  function defined by $Y_j/Y_i$ on
$\tilde{C}_i$. We set $v_{\tau,i} = \tau y_{j,i}+y_{k,i}$, where
 $\{i,j,k\} = \{1,2,3\}$, $j <k$ and $\tau \in {\mbox{\Bdd Z}}$.
Further, we denote by $\psi_{i,j,k}$ the function defined on
$\tilde{C}_i$ by the fraction $\psi_{i,j}/\psi_{i,k}$. By
Proposition 1, there is a monic polynomial $g_i(T) \in K[T]$ such
that the function $\tilde{v}_{\tau,i} =
g_i(\psi_{i,j,i})v_{\tau,i}$ is integral over
$K[\psi_{i,j,i},\psi_{i,k,i}]$,  $\deg g_i \leq \bar{N}$ and
$$H(g_i)  < \Lambda_1(\rho,M,N,\bar{N}) H(F)^{\bar{N}} H(\tilde{F}_i)^{MN}
H(\Psi_i)^{N\bar{N}}.$$
 The zeros of $g_i(T)$ are the elements
$\psi_{i,j,i}(R)$, where $R \in {\psi_i}^{-1}(D_C(X_i))\cap
V_{\tilde{C}_i}(Y_i)$. Moreover, there is $\Pi_{\tau,
i}(X_j,X_k,U)\in K[X_j,X_k,U]$ such that
$\Pi_{\tau,i}(\psi_{i,j,i},\psi_{i,k,i},U)$ is the minimal
polynomial of $\tilde{v}_{\tau,i}$ over the ring
$K[\psi_{i,j,i},\psi_{i,k,i}]$. Write
$$\Pi_{\tau, i}(X_j,X_k,U) =  U^{\nu}+
\pi_{\tau,i,1}(X_j,X_k)U^{\nu-1}+{\cdots}+\pi_{\tau,i,\nu}(X_j,X_k).$$
We have $\nu  \leq m$, $\deg \pi_{\tau,i,l} < 11MN^4\bar{N}^2$ $(l
= 1,\ldots,\nu)$ and
$$H(\Pi_{\tau,i}) < \Lambda_8(\tau,m,M,N, \bar{N})
(H(F)^{6N^2\bar{N}} H(\Psi_i)^{\bar{N}}
H(\tilde{F}_i)^M)^{240mM^3N^{12}\bar{N}^5}.$$ By \cite[Lemma
4.4]{Draziotis},  $H(\tilde{F}_i) < \Lambda_9(\bar{N},s(i))
H(\bar{F})$ and $ H(\Psi_i) < \Lambda_{10}(M,s(i)) H(\Phi).$ It
follows that $H(g_i)$ and $H(\Pi_{\tau,i})$ satisfy inequalities
as above having $H(\bar{F})$ and $H(\Phi)$ in place of
$H(\tilde{F}_i)$ and $H(\Psi_i)$ respectively.

The points $(z_1:z_2:z_3) \in D_C(X_i)$ with $z_i = 1$ and
$g_i(z_j) = 0$ belong to ${\phi}(V_{\bar{C}}(X_k+s(i)X_j))$. On
the other hand,  $P \in B_i$ and $B_i{\cap}
\phi(V_{\bar{C}}(X_k+s(i)X_j) = \O$. Hence, $g_i(a_j/a_i) \neq 0$
and so $v_{\tau,i}$ is defined at $\tilde{Q}$ $(i = 1,2,3)$.

Let $\psi_i^{*} : \overline{K}(C) \rightarrow
\overline{K}(\tilde{C}_i)$ be the field homomorphism associated to
the morphism $\psi_i$. As previously, there is a set $\Delta_i
\subset \mbox{\Bdd Z}$ with $|\Delta_i| \leq m(m-1)+2\bar{N}$ such
that for every integer $\tau \not \in \Delta_i$ we have
$\overline{K}(\tilde{C}_i) =
\psi_i^{*}(\overline{K}(C))(\tilde{v}_{\tau,i})$ (so $\nu = m$)
and $K(v_{\tau,i}(\tilde{Q})) = K(\tilde{Q}) = K(Q)$.

Let $\Sigma_{\tau,i}(X_j,X_k)$ be the discriminant of
$\Pi_{\tau,i}(X_j,X_k,U)$  with respect to $U$. We have $\deg
\Sigma_{\tau,i} \leq (2m-1)11\bar{N}^2N^4M$. We denote by $\Xi_i$
the set of points $(z_1:z_2:z_3) \in D_C(X_i)$ with $z_i = 1$, $
D_{r(i),i}(z_j,z_k) = 0$ and  $ \Sigma_{\tau,i}(z_j,z_k) =  0$,
for every  $\tau \in {\mbox{\Bdd Z}}$. Suppose that
$(z_1:z_2:z_3)\in \Xi_i$ with $z_i = 1$. Then, for every $\tau \in
\mbox{\Bdd Z}$, $\Pi_{\tau,i}(z_j,z_k,U)$ has at most $m-1$
distinct roots. If $g_i(z_j) \neq 0$, then there are $m$ distinct
points $Q_t \in \tilde{\phi}_i^{-1}(z_1:z_2:z_3)$ $(t = 1,\ldots,
m)$ and $\tau_0 \in \mbox{\Bdd Z}$ such that
$\tilde{v}_{\tau_0,i}(Q_p) \neq \tilde{v}_{\tau_0,i}(Q_q)$ for $p
\neq q$.  Thus, $\Pi_{\tau_0,i}(z_j,z_k,U)$ has  $m$ distinct
roots which is a contradiction. Hence $g_i(z_j) =  0$. Then
$(z_1:z_2:z_3) \in {\phi}(V_{\bar{C}}(X_k+s(i)X_j) \cap B_i = \O$
which is a contradiction. So, for every $(z_j,z_k)\in
\overline{K}^2$ with $D_{r(i),i}(z_j,z_k) =F_i(z_j,z_k) = 0$, the
polynomial in $\tau$, $\Sigma_{\tau,i}(z_j,z_k)$, is not zero.

Since $\tilde{v}_{\tau,i}$ is a root of
$\Pi_{\tau,i}(\psi_{i,j,i},\psi_{i,k,i},U)$,
$\pi_{\tau,i,l}(\psi_{i,j,i},\psi_{i,k,i})$, as polynomial in
$\tau$, has degree $\leq l$. Hence, the degree in $\tau$ of
$\Sigma_{\tau,i}(\psi_{i,j,i},\psi_{i,k,i})$ is $\leq (2m-1)m$.
So, for every $(z_1,z_2,z_3) \in J_{r(i),i}$ with $z_i = 1$ there
are at most $(2m-1)m$ integers $\tau$, such that
$\Sigma_{\tau,i}(z_j,z_k) = 0$. Thus, there is  $\tau(i) \in
{\mbox{\Bdd Z}}$ with $|\tau(i)| < 22m^3M\bar{N}^2N^5$, such that
$\overline{K}(\tilde{C}_i) =
\psi_i^{*}(\overline{K}(C))(\tilde{v}_{\tau(i),i})$ (so $\nu =
m$), $K(v_{\tau(i),i}(\tilde{Q})) =  K(Q)$ and
 for every $(z_1,z_2,z_3) \in J_{r(i),i}$
with $z_i = 1$ we have $\Sigma_{\tau(i),i}(z_j,z_k) \neq  0$.

Let $D_{\rho,i}^1$ and $\Sigma_{\tau,i}^1$ be two points in the
projective space having as coordinates 1 and the coefficients  of
$D_{\rho,i}(X_j,X_k)$ and  $\Sigma_{\tau,i}(X_j,X_k)$,
respectively. By \cite[Lemma 4.2]{Draziotis}, we have
$$H(D_{\rho,i}^1) < m^{3m-1} (11MN^4\bar{N}^2)^{4m-2}
H(P_{\rho,i})^{2m-1},$$
$$H(\Sigma_{\tau,i}^1) < m^{3m-1}
(11MN^4\bar{N}^2)^{4m-2} H(\Pi_{\tau,i})^{2m-1}.$$
We may assume,
without loss of generality, that one of the coefficients of $F$ is
 1.  By \cite[Lemma 4.1]{Draziotis}, there are positive integers
$a_{\rho,i},b_{\rho,i},c$ with   $$c \leq H_K(F)^{2N^2}, \ \
a_{\rho,i}  \leq H_K(P_{\rho,i})^{61mM^2\bar{N}^4N^8}, \ \
b_{\rho,i}  \leq H_{K}(\Pi_{\tau,i})^{61mM^2\bar{N}^4N^8}$$ such
that $a_{\rho,i}P_{\rho,i}(X_j,X_k,U)$,
$b_{\rho,i}\Pi_{\rho,i}(X_j,X_k,U)$ and $cF_i(X_j,X_k)$ have all
theirs coefficients in $O_K$.
 So, $a_{\rho,i}^{2m-2}D_{\rho,i}(X_j,X_k),
b_{\rho,i}^{2m-2}{\Sigma}_{\tau,i}(X_j,X_k) \in O_{K}[X_j,X_k]$.
Since $D_{r(i),i}(X_j,X_k)$, $\Sigma_{\tau(i),i}(X_j,X_k)$ and
$F_i(X_j,X_k)$ have no common zero, \cite[Lemma 2.9]{Draziotis}
implies that there are $A_{i,s} \in O_K[X_j,X_k]$ $(s = 1,2,3)$
and $A_i \in O_K \setminus \{0\}$ such that
$$A_{i,1}a_{r(i),i}^{2m-1}D_{r(i),i}+A_{i,2}b_{\tau(i),i}^{2m-1}{\Sigma}_{\tau(i),i}+A_{i,3}cF_i = A_i.$$
Furthermore, for every archimedean absolute value $|{\cdot}|_{v}$
of $K$ we have
$$|A_i|_{v} \leq ((\delta+1)(\delta+2)/2)! |E_i|_v^{(\delta+1)(\delta+2)/2},$$
where $\delta = 11MN^5\bar{N}^2$ and $E_i$ is a point of the
projective space with coordinates the coefficients of
$a_{r(i),i}^{2m-1}D_{r(i),i}$,
$b_{\tau(i),i}^{2m-1}{\Sigma}_{\tau(i),i}$ and $cF_i$. The bounds
for $a_{r(i),i}$, $b_{\tau(i),i}$, $c$, $H(D_{r(i),i}^1)$,
$H(\Sigma_{\tau(i),i}^1)$, $H(P_{r(i),i})$ and
$H({\Pi}_{\tau(i),i})$ give
$$|N_{K}(A_i)| <  \Lambda_{11}(d,m,M,N,\bar{N}) (H(F)^{6N^2\bar{N}} H(\Phi_i)^{\bar{N}}
H(\bar{F})^M)^{\lambda d m^3 M^7N^{30}\bar{N}^{13}},$$ where
$\lambda$ is a numerical constant.

Let $p_i = (a_j/a_i,a_k/a_i)$.  Since  $D_{r(i),i}(X_j,X_k)$,
$\Sigma_{\tau(i),i}(X_j,X_k)$ and $F_i(X_j,X_k)$ have no common
zero, we have either $D_{r(i),i}(p_i) \neq 0$ or
$\Sigma_{\tau(i),i}(p_i) \neq 0$. Let $S$ be the set of prime
ideals of $O_K$ dividing $A_1A_2A_3$.  Suppose that $\wp$ is a
prime ideal of $O_K$ with $\wp \not \in S$. Then there is $i \in
\{1,2,3\}$ such that $a_j/a_i,a_k/a_i \in O_{K,\wp}$. Put $L =
K(Q)$ and $\xi = [L:K]$. We have  $L = K(u_{r(i),i}(Q)) =
K(v_{\tau(i),i}(\tilde{Q}))$. We denote by $O_{K,\wp}$ the local
ring at $\wp$, by $\tilde{\wp}$ the prime ideal of $O_{K,\wp}$
generated by $\wp$ and by $D_{\wp}$ the discriminant of the
integral closure of $O_{K,\wp}$ into $L$ over $O_{K,\wp}$. Since
$\wp$ does not divide $A_i$, it follows that either
$a_{r(i),i}^{2m-1}D_{r(i),i}(p_i)$ or
$b_{\tau(i),i}^{2m-1}{\Sigma}_{\tau(i),i}(p_i)$ is not divisible
by $\tilde{\wp}$  (into $O_{K,\wp}$). If $\tilde{\wp}$ does not
divide $a_{r(i),i}^{2m-1}D_{r(i),i}(p_i)$, then
 $\tilde{\wp}$ does not divide  $a_{r(i),i}$ and $a_{r(i),i}^{2m-2}D_{r(i),i}(p_i)$.
Thus $a_{r(i),i}$ is a unit in $O_{K,\wp}$ and so $u =
u_{r(i),i}(Q)$ is integral over $O_{K,\wp}$. Then $D_{\wp}$
divides the discriminant $D(1,u, \ldots ,u^{\xi-1})$ of $1,u,
\ldots,u^{\xi-1}$ into  $O_{K,\wp}$. Further,
$D(1,u,\ldots,u^{\xi-1})$ divides
$a_{r(i),i}^{2m-2}D_{r(i),i}(p_i)$. Since $\tilde{\wp}$ does not
divide $a_{r(i),i}^{2m-2}D_{r(i),i}(p_i)$, $\tilde{\wp}$ does not
divide $D_{\wp}$. Thus, $\wp$ is not ramified into $L$. If
$\tilde{\wp}$ does not divide
$b_{\tau(i),i}^{2m-1}{\Sigma}_{\tau(i),i}(p_i)$, then we have the
same result.  By \cite[Lemma 4.3]{Draziotis},
$$N_K(D_{L/K}) <
 \prod_{\wp{\in}S}N_K(\wp)^{m-1} \exp (2m^2d) \leq N_K(A_1A_2A_3)^{m-1} \exp (2m^{2}d).$$
Using the estimates for $ N_K(A_i)$, the result follows.

\begin{flushleft}
Dimitrios Poulakis, \\
Aristotle University of Thessaloniki, \\
Department of Mathematics,\\
54124 Thessaloniki, Greece\\
Email Adress: poulakis@math.auth.gr\\

\   \\

Konstantinos Draziotis\\
Kromnis 33\\
54454 Thessaloniki, Greece\\
Email Adress: drazioti@math.auth.gr
\end{flushleft}


\begin{thebibliography}{99}

\bibitem{Bilu} Y. Bilu, Effective Analysis of Integral Points on Algebraic Curves, Ph. D.
Thesis, Beer Sheva, 1993.

\bibitem{Chevalley} C. Chevalley, Un th\'eor\`eme d'arithm\'etique sur les courbes alg\'ebriques,
{\it C. R. Acad. Sci. Paris} 195 (1932), 570-572.

\bibitem{Draziotis} K. Draziotis and D. Poulakis,  An Explicit Chevalley-Weil Theorem  for Affine Plane Curves,
Rocky Mountain Journal of Mathematics, Rocky Mountain Journal of
Mathematics, 39(1) (2009), 49-70.

\bibitem{Hindry} M. Hindry - J. Silverman, {\it Diophantine Geometry}, New-York Inc.: Springer-Verlag 2000.

\bibitem{Lang 83} S. Lang, {\it Diophantine Geometry,} Springer Verlag 1983.

\bibitem{Poulakis1} D. Poulakis, Polynomial bounds for the solutions of a class of Diophantine equations,
 {\it J. Number Theory}, 66, 2 (1997), 271-281.

\bibitem{Poulakis2} D. Poulakis, Integer points on algebraic curves with exceptional units, {\it J. Austral. Math. Soc.
(Series A)} 63 (1997), 145-164.

\bibitem{Shafarevich} I. Shafarevich, {\it Basic Algebraic Geometry}, Berlin-Heidelberg-New York: Springer Verlag 1977.

\bibitem{Silverman} J. H. Silverman, {\it The Arithmetic of Elliptic Curves,} Springer Verlag 1986.

\bibitem{Weil} A. Weil, Arithm\'etique et g\'eom\'etrie sur les vari\'et\'es alg\'ebriques, {\it Act. Sci. et Ind.} No 206, Paris:
Hermann 1935.

\end{thebibliography}
\end{document}